\documentclass[journal]{IEEEtran}
\usepackage{cite}
\usepackage{float}
\usepackage{graphicx} 
\usepackage{tikz}
\usepackage{pgfplots}
\usepackage{color}
\usepackage{xcolor}
\usepackage{amstext,epsfig,amssymb,amsbsy,amsmath}
\usepackage{mathtools}
\usepackage[output-decimal-marker={.}]{siunitx}
\usepackage{array}
\usepackage{tabularx}
\usepackage{multirow}
\usepackage{multicol}
\usepackage{setspace}
\usepackage[caption=false,font=footnotesize]{subfig}
\usepackage{wrapfig}
\usepackage{url}
\usepackage{hyperref}
\usepackage{cleveref}
\usepackage{balance}

\graphicspath{{Figure/}}

\expandafter\def\expandafter\UrlBreaks\expandafter{\UrlBreaks
  \do\a\do\b\do\c\do\d\do\e\do\f\do\g\do\h\do\i\do\j%
  \do\k\do\l\do\m\do\n\do\o\do\p\do\q\do\r\do\s\do\t%
  \do\u\do\v\do\w\do\x\do\y\do\z\do\A\do\B\do\C\do\D%
  \do\E\do\F\do\G\do\H\do\I\do\J\do\K\do\L\do\M\do\N%
  \do\O\do\P\do\Q\do\R\do\S\do\T\do\U\do\V\do\W\do\X%
  \do\Y\do\Z}

\newcommand\Tstrut{\rule{0pt}{2.6ex}}         
\newcommand\Bstrut{\rule[-0.9ex]{0pt}{0pt}}   

\definecolor{gray1}{rgb}{0.23529,0.23529,0.23529}%
\definecolor{gray2}{rgb}{0.49412,0.49412,0.49412}%
\definecolor{gray3}{rgb}{0.86275,0.86275,0.86275}%
\definecolor{orange1}{rgb}{0.87059,0.49020,0.00000}%
\definecolor{blue1}{rgb}{0.20392,0.30196,0.49412}%
\definecolor{blue2}{rgb}{0.72941,0.83137,0.95686}%
\definecolor{blue3}{rgb}{0.15294,0.22745,0.37255}%
\definecolor{orange2}{rgb}{1.00000,0.60000,0.00000}%
\definecolor{orange3}{rgb}{1.00000,0.89020,0.66667}%

\pgfplotsset{
    layers/my layer set/.define layer set={
        background,
        main,
        foreground
    }{},
    set layers=my layer set,
}

\begin{document}

\bstctlcite{IEEEexample:BSTcontrol}

\title{Market Integration of HVDC Lines: internalizing HVDC losses in market clearing}%

\author{Andrea~Tosatto,~\IEEEmembership{Student~Member,~IEEE,}
        Tilman~Weckesser,~\IEEEmembership{Member,~IEEE,}
        and~Spyros~Chatzivasileiadis,~\IEEEmembership{Senior~Member,~IEEE}%
\thanks{A. Tosatto and S. Chatzivasileiadis are with the Technical University of Denmark, Department of Electrical Engineering, Kgs. Lyngby, Denmark (emails: \{antosat,spchatz\}@elektro.dtu.dk). T. Weckesser is with Dansk Energi, Frederiksberg C, Denmark (email: twe@danskenergi.dk). This work was supported by the multiDC project, funded by Innovation Fund Denmark under Grant 6154-00020B. Digital Object Identifier 10.1109/TPWRS.2019.2932184.}}%

\maketitle

\vspace{-1em}
\begin{abstract}
Moving towards regional Supergrids, an increasing number of interconnections are formed by High Voltage Direct Current (HVDC) lines. Currently, in most regions, HVDC losses are not considered in market operations, resulting in additional costs for Transmission System Operators (TSOs). Nordic TSOs have proposed the introduction of HVDC loss factors in the market clearing algorithm, to account for the cost of losses and avoid HVDC flows between zones with zero price difference. In this paper, we introduce a rigorous framework to assess the introduction of loss factors, in particular HVDC loss factors, in nodal and zonal pricing markets. First, we focus on the identification of an appropriate loss factor. We propose and compare three different models: constant, linear, and piecewise linear. Second, we introduce formulations to include losses in market clearing algorithms. Carrying numerical tests for a whole year, we find that accounting only for HVDC or AC losses may lead to lower social welfare for a non-negligible amount of time. To counter this, this paper introduces a framework for including both AC and HVDC losses in a zonal or nodal pricing environment. We show both theoretically and through simulations that such a framework is guaranteed to increase social welfare.
\end{abstract}


\vspace{-0.5em}
\section{Introduction}\label{sec:1}
\IEEEPARstart{W}{ith} the progress made in the field of power electronics in the past decades, High Voltage Direct Current (HVDC) lines are now considered an attractive alternative to AC lines. Indeed, compared to AC, the transmission of power in the DC form presents several benefits, such as lower power losses beyond a certain distance, possibility of connecting asynchronous areas, full controllability of the power flows, no need of reactive power compensation, and others \cite{2_1, 2_2, 2_3}. All these features make HVDC lines particularly convenient in those applications where bulk power has to be transmitted over long distances. Consequently, contrary to AC interconnections, which usually span only a few hundred meters, HVDC interconnectors are often hundreds of kilometers long. Thus, when considering the operation of such long HVDC lines, the cost of thermal losses becomes non-negligible and the question that arises is: who should bear these costs?

\input{fig_HVDCmap.tex}

Ideally, the operation costs of transmission systems should be covered by generating companies and consumers through the market mechanisms. In the US, for example, CAISO and PJM Interconnection LLC include in the energy price the marginal cost of losses using Marginal Loss Sensitivity Factors (MLFs), which show the marginal increase in system losses due to a marginal increase in power injections at a specific location \cite{2_5, 2_6}. In Australia, losses on interconnectors are calculated using inter-zonal loss factors and, in a similar way, market participants within a bidding zone are charged based on intra-zonal loss factors \cite{2_7}, considering losses occurring between the Regional Reference Node and their point of connection.

In Europe, the Price Coupling of Regions (PCR) project aims at coupling twenty-five different countries whose electricity markets are operated by eight Power Exchanges \cite{2_8}. To avoid excessive complexity, market operators use a simplified model \cite{2_9} for determining power exchanges between different regions. Although losses on interconnectors could be considered, for the majority this is not done. As a result, the revenues of the Transmission System Operators (TSOs) through the market are insufficient to cover the extra costs of losses; grid tariffs are introduced to fill this gap, among others.

Different TSOs follow different practices in order to include the cost of losses in the grid tariffs. In certain regions, once the market has been cleared, an ex-post settlement is reached and the cost of losses is allocated across generators and loads, using sensitivity and transmission loss factors \cite{2_10}. Other TSOs estimate the total losses with an ex-ante calculation using offline models and the cost of losses are included in the grid tariffs. The share of losses among generators and loads varies from country to country, and usually loads carry a higher share to allow generators to be more competitive in the European Market \cite{2_11}. In addition to internal losses, losses on interconnectors are handled through special agreements between TSOs: losses are usually shared equally and each TSO bids in the day-ahead market for its share of losses. Concerning transit flows, the Inter-TSO Compensation Mechanism (ITC) aims at compensating the use of infrastructures and losses caused by hosting transit flows \cite{2_12}. However, it is specified that ``There should be no specific network charge for individual transactions for declared transit of energy'' \cite{2_13}, meaning that HVDC flows do not fit in this mechanism.

The problem of losses arises, especially for HVDC lines, when the price differences among zones are very small. This happens often in Scandinavia. For example, the price difference between Denmark (DK1) and Norway (NO2) has been zero for more than 4000 hours during 2018 \cite{2_14}. These two areas are connected by a 240-km long HVDC line, Skagerrak. If we consider the power exchanges during these hours, losses have cost more than 4 million Euros in 2018 while no revenue has been obtained from the electricity trade. Considering the large number of HVDC connections in Europe that face a similar situation, and the increasing number of new projects (\figurename~\ref{fig:1_map}), the cost of losses amounts to tens of millions of euros.

To deal with this problem, some TSOs are considering to internalize the cost of losses in the market clearing procedure, moving from an ``explicit'' to an ``implicit grid loss'' calculation. In \cite{2_15}, the TSOs of the Nordic Capacity Calculation Region (CCR) propose to include loss factors for only HVDC interconnectors in the market clearing, as HVDC losses are substantial and HVDC flows are fully-controllable. Through that, power flows among zones would only be allowed if the price difference is greater than the marginal cost of losses. In \cite{2_15}, the Nordic TSOs present the results of different simulations with implicit grid losses implemented on some of the HVDC interconnectors in the Nordic area. The following question arises: is the introduction of loss factors for only HVDC interconnectors the best possible action?

The aim of this paper is to introduce a rigorous framework for analyzing the inclusion of losses in the market clearing. More specifically, the contributions of this paper are:
\begin{itemize}
    \item the introduction of a framework to assess how incorporating the losses of AC and HVDC lines in market clearing affects the market outcome;
    \item the investigation of different loss factor formulations and their impact, while maintaining the linear formulation of the market clearing algorithm;
    \item a detailed method on how to include cross-border AC losses and intra-zonal losses in a zonal market;
    \item an analytical proof on how the inclusion of AC and/or HVDC losses impacts the market clearing outcome.
\end{itemize}

The paper is organized as follows. \mbox{Section \ref{sec:2}} outlines the modeling of HVDC interconnectors and AC grids for the Optimal Power Flow (OPF), and  \mbox{Section \ref{sec:3}} describes the market clearing algorithm for nodal and zonal markets. In \mbox{Section \ref{sec:4}}, we propose different loss factor formulations for HVDC and AC lines and analyze their properties. In \mbox{Section \ref{sec:5}}, we propose a methodology to derive loss factors for regional markets, considering also losses due to cross-border flows. \mbox{Section \ref{sec:6}} presents numerical results on a 4-area 96-bus system and \mbox{Section \ref{sec:7}} concludes.
\vspace{-0.5em}
\section{Transmission line modeling}\label{sec:2}
Due to the non-linear nature of power flow equations, most of the market clearing algorithms use a simplified model of the power system, following a ``DC power flow'' approximation \cite{3_1}: line resistances are assumed considerably smaller than line reactances, thus the transmission network is modeled using only the imaginary part of line impedances and no active power losses are implicitly calculated. Moreover, voltage magnitudes are assumed close to 1 p.u., thus line flows are determined only by the angle differences between nodes. In the following, the simplified model for HVDC and AC lines is presented.

\vspace{-1em}
\subsection{Point-to-point HVDC connections}\label{ssec:2_1}

\input{fig_VSCsimpl.tex}

An HVDC point-to-point connection consists of two converters connected through a DC power cable. The two converters are connected to AC systems, and the way they are modelled depends on the technology used for the conversion, that is Line-Commutated Converter (LCC) or Voltage-Source Converter (VSC).

Under the aforementioned assumptions, the complete model of HVDC lines (that can be found in \cite{3_2}) can be simplified, as shown in \figurename~\ref{fig:2_VSCsimplified}. In this model, all components inside the converter stations are substituted with an AC voltage source and the DC system is not included. With these modifications, the model is lossless and the power flowing over the line is equal to the power sent and received at the connected nodes. 

If we indicate with $f^{\textsc{dc}}_l$ the power flowing over line $l$, the power balance equation becomes:
\begin{equation}\label{eq:2_1_1}
    \sum_{n} \text{I}^{\textsc{dc}}_{n,l} \cdot f^{\textsc{dc}}_l = 0,
\end{equation}

where $\text{I}^{\textsc{dc}}$ is the HVDC line incidence matrix, defined as:
\begin{equation}\label{eq:2_1_2}
    \text{I}^{\textsc{dc}}_{n,l} = 
    \begin{cases}
    1, & \text{if bus $n$ is the receiving bus of line $l$} \\
    -1, & \text{if bus $n$ is the sending bus of line $l$} \\
    0, & \text{otherwise.}
    \end{cases}
\end{equation}
$\text{I}^{\textsc{dc}}$ is a $n_{bus} \times n_{line\textsc{dc}}$ matrix, where $n_{line\textsc{dc}}$ is the number of HVDC lines in the system and $n_{bus}$ the number of nodes. 

This HVDC model is a simplified version and is used only for the determination of power exchanges between bidding zones during the market clearing. The complete HVDC model, as outlined in \cite{3_2}, is used for offline calculation of losses, available transmission capacity and for security considerations.

\vspace{-1em}
\subsection{Meshed AC grids}\label{ssec:2_2}

AC lines are generally modeled with a $\pi$-equivalent model, consisting of an electrical impedance ($R+jX$) between the sending and receiving nodes and two shunt capacitances ($j\frac{b_{sh}}{2}$) at the connected nodes \cite{3_3}.

As mentioned above, by neglecting line resistances and shunt elements, AC lines can be modeled by their line susceptances, resulting in the simplified model used for DC power flow studies \cite{3_1}.

Contrary to point-to-point connections, in a meshed grid line flows are not free decision variables, but functions of the power injected at each node:
\begin{equation}\label{eq:2_2_1}
    \boldsymbol{f}^{\textsc{ac}} = \textbf{PTDF} \cdot \boldsymbol{P}^{\textsc{inj}}
\end{equation}

The Power Transfer Distribution Factor (PTDF) matrix shows the marginal variation in the power flows due to a marginal variation in the power injections and it is used to calculate how power flows are distributed among transmission lines.

For a power system with $n_{line}$ lines and $n_{bus}$ buses, the PTDF matrix is an $n_{line}\times n_{bus}$ matrix and can be calculated as:
\begin{equation}\label{eq:2_2_2}
    \textbf{PTDF} = \boldsymbol{B}_{line} \boldsymbol{\tilde{B}}_{bus}^{-1}
\end{equation}
where $\boldsymbol{B}_{line}$ is the line susceptance matrix and $\boldsymbol{\tilde{B}}_{bus}^{-1}$ is the inverse of the bus susceptance matrix after removing the row and the column corresponding to the slack bus \cite{3_1}. 
\section{Market clearing algorithm}\label{sec:3}
Under the assumption of perfectly competitive electricity markets, the market-clearing outcome is a \textit{Nash equilibrium}, that is a state in which none of the producers or consumers can increase its profit by deviating from the equilibrium, i.e. changing unilaterally its schedule. The equilibrium model of electricity markets consists of four blocks, each one representing a different market participant. In the first block, each producer maximizes its profit from the sale of energy. Similarly, in the second block, each elastic load maximizes its profit from the purchase of energy. The third block represents the profit-maximization problem of the transmission system operator, who seeks to maximize the profit from the trade of electricity among different areas. Finally, the last block consists of the common market constraints, i.e. power balance equations. The formulation of an optimization problem for each market participant gives the freedom to arbitrarily change their objective functions and, thus, include losses. 

All the optimization problems in the equilibrium model are linear and convex, thus it is possible to substitute them with their Karush-Kuhn-Tucker (KKT) optimality conditions. By doing so, the equilibrium model is recast as a mixed-complementarity problem (MCP) including the KKT conditions and the linking constraints. MCPs can be solved using the PATH solver on GAMS, or other similar solvers. Another possibility, under certain circumstances, is to recast the MCP as a single optimization problem. This is possible only when there exists an optimization problem with the same optimality conditions as the original MCP \cite{4_1}. In this case, since all the market participants are price-takers, the dual variables (that influences the market prices) are parameters in their optimization problems. Thus, it is possible to recast the MCP problem as a single optimization problem \cite{4_1}.

In order to obtain a feasible dispatch, i.e. a set of injections that does not violate any network constraint, the transmission network is included in the market model. Depending on how the network is modeled, different pricing mechanisms are possible. In the following, a brief description of nodal and zonal pricing markets is given.

\vspace{-1em}
\subsection{Nodal pricing markets}\label{ssec:3_1}
In a nodal pricing system, all transmission lines and transformer substations are included in the network model. Locational Marginal Prices (LMPs) are defined for each node of the network, and generators and loads are subjected to a different price according to the substation they are connected to \cite{4_2}. This is the case of different markets in the US, e.g. the Californian electricity market operated by CAISO (California Independent System Operator) or the market operated by PJM Interconnection LLC (Pennsylvania-NewJersey-Maryland) \cite{2_5, 4_3}, and other markets, such as New Zealand's Exchange (NZX) \cite{4_4}.

\onehalfspacing
In its simplest form, the market-clearing problem can be formulated as the following optimization problem:
\vspace{0.4em}
\small
\begin{subequations}\label{eq:3_1_1}
    \begin{align}
        \underset{\boldsymbol{g}, \boldsymbol{d}, \boldsymbol{f}^{\textsc{dc}}}{\text{max}} & \quad \boldsymbol{u}^\intercal \boldsymbol{d} - \boldsymbol{c}^\intercal \boldsymbol{g} \label{3_1:obj} \\
        \text{s.t.} & \quad \underline{\boldsymbol{G}} \leq \boldsymbol{g} \leq \overline{\boldsymbol{G}} \label{3_1:gen} \\
        & \quad \underline{\boldsymbol{D}} \leq \boldsymbol{d} \leq \overline{\boldsymbol{D}} &&  \label{3_1:con} \\
        & \quad \boldsymbol{f}^{\textsc{ac}} = \textbf{PTDF} \cdot ( \boldsymbol{\text{I}}^{\textsc{g}}\boldsymbol{g} - \boldsymbol{\text{I}}^{\textsc{d}}\boldsymbol{d} - \boldsymbol{\widetilde{p}}^{\,loss\textsc{n}} + \boldsymbol{\text{I}}^{\textsc{dc}}\boldsymbol{f}^{\textsc{dc}} ) \label{3_1:flowAC} \\
        & \quad -\overline{\boldsymbol{F}}^{\textsc{ac}} \leq \boldsymbol{f}^{\textsc{ac}} \leq \overline{\boldsymbol{F}}^{\textsc{ac}} \,\,\, : \,\,\, \boldsymbol{\underline{\mu}}^{\textsc{ac}}, \boldsymbol{\overline{\mu}}^{\,\textsc{ac}} \label{3_1:limitAC} \\
        & \quad -\overline{\boldsymbol{F}}^{\textsc{dc}} \leq \boldsymbol{f}^{\textsc{dc}} \leq \overline{\boldsymbol{F}}^{\textsc{dc}} \label{3_1:limitDC} \\
        & \quad \sum_{j} d_{j} + \sum_{n} \widetilde{p}^{\,loss\textsc{n}}_n - \sum_{i} g_{i} = 0 \,\,\, : \,\,\, \lambda \label{3_1:balance}
    \end{align}
\end{subequations} %
\vspace{-1.4em}
\normalsize \singlespacing \noindent %
where $\boldsymbol{u}$ and $\boldsymbol{c}$ are respectively load utilities and generator costs, $\boldsymbol{g}$ and $\boldsymbol{d}$ are the output levels of generators and consumption of loads, $\boldsymbol{\text{I}}^{\textsc{g}}$ and $\boldsymbol{\text{I}}^{\textsc{d}}$ are the incidence matrices of generators and load, $\underline{\boldsymbol{G}}$, $\overline{\boldsymbol{G}}$, $\underline{\boldsymbol{D}}$ and $\overline{\boldsymbol{D}}$ are respectively the minimum and maximum generation and consumption of each generator and load, $\boldsymbol{f}_{\textsc{ac}}$ and $\boldsymbol{f}_{\textsc{dc}}$ are the power flows over AC and HVDC lines, $\overline{\boldsymbol{F}}^{\textsc{ac}}$ and $\overline{\boldsymbol{F}}^{\textsc{dc}}$ are the capacities of AC and HVDC lines, $\boldsymbol{\underline{\mu}}^{\textsc{ac}}$ and $\boldsymbol{\overline{\mu}}^{\,\textsc{ac}}$ are the lagrangian multipliers associated with AC line limits, $\lambda$ is the lagrangian multiplier associated with the power balance equation and $\boldsymbol{\widetilde{p}}^{\,loss\textsc{n}}$ are the nodal losses. For now, it is assumed that losses are just parameters in the optimization problem, which are estimated using off-line models before the market is cleared.

The objective of the market operator is to maximize the social welfare, expressed in \eqref{3_1:obj} as the difference between load pay-offs and generator costs. Constraints \eqref{3_1:gen} and \eqref{3_1:con} enforce the lower and the upper limits of generation and consumption, while constraints \eqref{3_1:limitAC} and \eqref{3_1:limitDC} ensure that line limits are not exceeded. The flows over AC interconnectors are defined through constraint \eqref{3_1:flowAC} using the PTDF matrix (see Section \ref{ssec:2_2}).

The LMPs are computed as follows \cite{4_5}:
\begin{equation}\label{eq:3_1_2}
    \text{LMP}_n = \lambda + \sum_l \text{PTDF}_{n,l} (\underline{\mu}^{\textsc{AC}}_{\,l} - \overline{\mu}^{\,\textsc{AC}}_l ) 
\end{equation}

\vspace{-1em}
\subsection{Zonal pricing markets}\label{ssec:3_2}

In a zonal pricing system, the network is split into price-zones in case of congestion on certain flowgates. The intra-zonal network is not included in the model, and a single price per zone is defined. The main difference between nodal and zonal pricing is that, in case of congestion, in a nodal pricing market all the nodes are subjected to different prices, while in a zonal pricing market price differences arise only among zones, with all generators and loads subjected to their zonal price \cite{4_2}. An example of this pricing system is the Australian electricity market operated by AEMO (Australian Energy Market Operator) \cite{4_6}. An evolution of zonal pricing is the Flow-Based Market Coupling (FBMC), which aims at coupling different independent markets. FBMC includes two clearing processes: first the energy market clearing, where a clearing price per zone is determined according to the internal power exchanges, and second, the import and export trades via the interconnections \cite{4_2}. As for zonal pricing, the intra-zonal flows are not represented in the model; in addition, cross-border lines to another zone are aggregated into a single equivalent interconnector. This is the underlying concept of the Price Coupling of Regions (PCR) project of the European Power Exchanges (EPEX) \cite{4_7}.

As a result, in zonal pricing market, the PTDF matrix becomes an $n_{line}\times n_{zone}$ matrix and must be \emph{estimated} taking into consideration the intra-zonal networks that are omitted in the market model, so that the resulting flows from the market clearing can resemble the actual ones. The estimation of PTDFs is based on statistical factors related to flows on the bidding zone borders under different load and generation conditions \cite{4_8}. The PTDF matrix is calculated as follows. One at a time, the output of all generators is increased by \SI{1}{\mega\watt}: for each generator, power flow analyses are carried out considering the extra megawatt consumed at a different bus every time. For all the generation patterns and load conditions, the marginal variation of the power flows on the interconnectors is calculated. At the end, the PTDFs are estimated by statistical analysis using linear regression.

Once the \emph{zonal} PTDF is calculated, the market is cleared solving the optimization problem \eqref{eq:3_1_1}. Zonal prices are still calculated with Eq. \eqref{eq:3_1_2}, but they refer to regions and not to single buses.
\vspace{-0.5em}
\section{Loss factor formulations}\label{sec:4}
\subsection{HVDC losses}\label{ssec:4_1}

The power losses of an HVDC link can be calculated as the sum of the losses in the two converter stations plus the losses on the DC cable. The latter are the ohmic losses due to the resistance of the cable, calculated as: 
\begin{equation}\label{eq:4_1_1}
    p^{cable} = R\,\left|I^{line}\right|^2,
\end{equation}
where $R$ is the resistance of the cable and $\left|I^{line}\right|$ is the magnitude of the line current.

In an HVDC converter station, most of the losses are due to the transformer, the AC filter, the phase reactor and the converter. Transformer losses are calculated in the same way as for conventional power transformers in AC grids, and they can be divided into iron losses (no-load losses) and copper losses (load losses). Due to the high harmonic content of the current, however, losses tent to be higher than in conventional transformers \cite{5_1}. Losses in the filters and in the phase reactor are due to their parasitic resistances which are modeled as equivalent series resistances. For load flow analysis, the above mentioned losses are calculated by including the impedances of these elements in the admittance matrix of the system.

Converter losses can be divided into switching and conduction losses: the first are caused by the turn-on and turn-off of power electronic devices, the second are the ohmic losses caused by their parasitic resistances during the on-state mode. Power losses can vary significantly depending on the converter technology: LCC converters use thyristors as switching devices, while VSCs use IGBT's. Thyristors are semiconductor devices that can only be turned on by control action, resulting in a single commutation per cycle \cite{5_2}. With IGBT's, both turn-on and turn-off can be controlled, giving an additional degree of freedom. This controllability comes with a price, since IGBT's are switched on and off many times (typically between 20 and 40) per cycle. For this reason, switching losses are significantly higher in VSC converters \cite{5_3}. With increasing number of IGBT's per arm, switching losses tend to decrease. In a three-level topology, for example, the number of commutations per cycle is half compared to a two-level topology. With Modular Multi-level Converters (MMCs), each valve is composed by several independent converter submodules containing two IGBT's connected in series. In each submodule, IGBT's are switched on and off only once per cycle, resulting in less switching losses \cite{5_4}. As a rule of thumb, one can say that a typical LCC-HVDC converter station has power losses of around 0.7\%, a VSC-HVDC converter station of around 2-3\% and MMC-HVDC converter station of around 1\% \cite{5_4}. The switching frequency does not only influence the converter losses, but has also an impact on the losses produced by the other devices. Indeed, switching losses increase with the switching frequency, but the harmonic distortion decreases. Since losses in other components depends on the RMS value of the current, their losses increase with its harmonic content, thus with lower switching frequencies. Concerning the operating mode, when a VSC converter is operating as a rectifier, diodes are conducting more frequently than IGBT's, and since IGBT's have higher conduction losses than diodes, the losses are lower compared to the inverter mode, when IGBT's are used more frequently \cite{5_1}. This does not happen with LCC converters, since the only difference between the two operating modes is the firing angle and there are always two valves conducting at a time.

The operation of the converter station requires also a certain number of auxiliary services, such as auxiliary power supply, valve cooling, air conditioning, fire protection, etc. According to \cite{5_6}, the total power consumption of auxiliary devices is around 0.1\% of the converter station rating. The losses of the remaining equipment, such as switchgear, instrument transformers, surge arresters, etc., are negligible compared to the above mentioned losses \cite{5_7}.

The modeling of losses through the calculation of equivalent impedances would require a detailed knowledge of all the individual loss contributions of these devices, thus, the converter station losses are commonly represented with the generalized loss model \cite{3_2}:
\begin{equation}\label{eq:4_1_2}
    p^{conv} = a\,\left|I^{conv}\right|^2 + b\,\left|I^{conv}\right| + c,
\end{equation}
where $a$, $b$ and $c$ are numerical parameters reflecting the quadratic, linear and constant dependence of the losses on the line current and $\left|I^{conv}\right|$ is the magnitude of the current flowing through the converter. The constant parameter represents the amount of losses that is produced also when the HVDC link is not operated, that is when the converter station is energized but the valves are blocked. With the right choice of $a$, $b$ and $c$, this generalized loss model is suitable for both VSC- and LCC-based HVDC converter stations and for different converter topologies. Examples of these parameters can be found in \cite{3_2, 5_1}.

The losses on an HVDC line are thus quadratic function of the current, as shown in \eqref{eq:4_1_1} and \eqref{eq:4_1_2}. However, as explained in Section \ref{sec:2}, many market-clearing algorithms use a simplified model that considers linear functions of active power. The first step towards a linear approximation of HVDC losses is to replace the line current in \eqref{eq:4_1_1} and the converter current in \eqref{eq:4_1_2} with the HVDC active power flow. Since no shunt elements are considered in the simplified model, then $\left|I^{conv}\right|=\left|I^{line}\right|=\left|I^{\textsc{dc}}\right|$. Working in the \textit{per unit} system, and assuming $|V|=1$~p.u. at each bus (which is the standard DC power flow approximation), then $\left|f^{\textsc{dc}}\right|=\left|I^{\textsc{dc}}\right|$. As a result, for the HVDC line $l$, the total losses can be approximated to:
\begin{equation}\label{eq:4_1_3}
    p^{loss\textsc{dc}}_l = A_l \left|f^{\textsc{dc}}_l\right|^2 + B_l \left|f^{\textsc{dc}}_l\right| + C_l.
\end{equation}
with $A_l = a^{inv}_{l}+a^{rect}_{l}+R_l$, $B_l = 2b_l$ and $C_l=2c_l$, where $R_l$ is the resistance of the cable and $a^{inv}_{l}$, $a^{rect}_{l}$, $b_l$ and $c_l$ are respectively the quadratic, linear and constant loss coefficients of the converter stations (one operating in inverter mode, the other in rectifier mode). Different linearization techniques are presented in Section \ref{ssec:4_3}.

\vspace{-1em}
\subsection{AC losses}\label{ssec:4_2}
Losses in AC grids are produced by a large number of devices; however, under the assumption of DC power flows, it is common practice to only include transmission lines and transformers in the network model.

Since no reactive power flows are considered, losses on a transmission line can be expressed as \cite{5_8, 5_10, 5_11}:
\begin{equation}\label{eq:4_2_1}
    p^{line} = R \left|I^{line}\right|^2
\end{equation}
where $R$ is the resistance of the line and $\left|I^{line}\right|$ the magnitude of the line current.

As described in Section \ref{ssec:4_1}, transformer losses can be divided into iron and copper losses, and are modeled through equivalent resistances:
\begin{equation}\label{eq:4_2_2}
    p^{tran} = R^{eq} \left|I^{tran}\right|^2
\end{equation}
where $R^{eq}$ is the equivalent resistance of the transformer and $\left|I^{tran}\right|$ is the current flowing through the transformer. However, when solving a DC power flow problem, no distinction is made between transformers and transmission lines. Moreover, using the per-unit system, the line or transformer current can be substituted by the active power flow. It is possible, thus, to express the losses occurring between two AC buses with the general loss function:
\begin{equation}\label{eq:4_2_3}
    p^{loss\textsc{ac}}_l = R_l \left|f^{\textsc{ac}}_l\right|^2
\end{equation}

\vspace{-1em}
\subsection{Linearization techniques}\label{ssec:4_3}
To avoid excessive complexity, most of the market clearing software (e.g. \cite{2_9} or \cite{5_19}) don't allow polynomial constraints with degree above 1. In this section, three linearization techniques are introduced: constant, linear and piecewise linear. 

\subsubsection{Constant loss factors}\label{sssec:4_3_1}

One possibility is to consider the losses constant:
\begin{equation}\label{eq:4_3_1}
    p^{loss}_l = \beta_l.
\end{equation}
The coefficient $\beta_l$ can be estimated considering losses during the maximum power flowing through the line, or losses occurring with a certain power flow. In the second case, the average power flowing on the line can be calculated considering a time window of one year. 

\subsubsection{Linear loss factors}\label{sssec:4_3_2}

If we consider linear dependence of losses on the power flow, the loss equation becomes:
\begin{equation}\label{eq:4_3_2}
    p^{loss}_l = \alpha_l \left|f_l\right| + \beta_l.
\end{equation}
Parameters $\alpha_l$ and $\beta_l$ can be estimated in different ways, e.g. using the least squares approach, connecting stand-by losses to maximum losses, linearizing around a certain range of flows, through the derivative at a certain flow, etc. 

\subsubsection{Piecewise linear loss factors}\label{sssec:4_3_3}

\onehalfspacing
A better approximation of losses is obtained by constructing a piecewise linear function. With $K$ segments, the loss equation becomes:
\small
\begin{equation}\label{eq:4_3_3}
    p^{loss}_l = 
    \begin{cases}
    \alpha_{1,l} \left|f_l\right| + \beta_{1,l}, & \text{if $\left|f_l\right|\leq f^*_1$} \\
    \qquad \quad \vdots & \\
    \alpha_{K,l} \left|f_l\right| + \beta_{K,l}, & \text{if $f^*_{K-1}\leq \left|f_l\right| \leq f^*_K$} \\
    \end{cases}
\end{equation}
\vspace{-2em}
\normalsize \singlespacing
\noindent
For each line segment $k$, parameters $\alpha_{k,l}$ and $\beta_{k,l}$ can be calculated in a similar way as explained above for linear loss factors.

\vspace{-1em}
\subsection{Inclusion of losses in the market clearing}\label{ssec:4_4}

Although convex, the absolute value operator is non-linear. For this reason, when added to problem \eqref{eq:3_1_1}, equation \eqref{eq:4_3_2} or each equality of \eqref{eq:4_3_3} is recast as two inequalities in the form of:
\begin{equation}\label{eq:4_4_1}
    \begin{aligned}
        & p^{loss}_l \geq \alpha_l f_l + \beta_l \,\,:\,\,\sigma^+_l & \,\,\, \forall l  \\
        & p^{loss}_l \geq \alpha_l (-f_l) + \beta_l \,\,:\,\,\sigma^-_l & \,\,\, \forall l   
    \end{aligned}
\end{equation}

Once losses are calculated, they are considered as an additional load and equally split between the buses at the sending and the receiving end. For this purpose, a \textit{loss distribution matrix} is defined as follows:
\begin{equation}\label{eq:4_4_2}
    D_{n,l} = 
    \begin{cases}
    0.5, & \text{if line $l$ is connected to bus $n$} \\
    0, & \text{otherwise}
    \end{cases}
\end{equation}
Nodal losses $\boldsymbol{p}^{loss\textsc{n}}$ and zonal losses $\boldsymbol{p}^{loss\textsc{z}}$ are now calculated as:
\begin{equation}\label{eq:4_4_3}
    \begin{aligned}
        & \boldsymbol{p}^{loss\textsc{n}} = \boldsymbol{D}^{\textsc{dc}} \cdot \boldsymbol{p}^{loss\textsc{dc}} + \boldsymbol{D}^{\textsc{ac}} \cdot \boldsymbol{p}^{loss\textsc{ac}} \\
        & \boldsymbol{p}^{loss\textsc{z}} = \boldsymbol{D}^{\textsc{dc}} \cdot \boldsymbol{p}^{loss\textsc{dc}} + \boldsymbol{D}^{\textsc{ac}} \cdot \boldsymbol{p}^{loss\textsc{ac}} + \boldsymbol{\widetilde{p}}^{\,intra}
    \end{aligned}
\end{equation}
where $\boldsymbol{p}^{loss\textsc{ac}}$ and $\boldsymbol{p}^{loss\textsc{dc}}$ are the losses on AC and HVDC lines, $\boldsymbol{D}^{\textsc{ac}}$ and $\boldsymbol{D}^{\textsc{dc}}$ are the AC and HVDC loss distribution matrices and $\boldsymbol{\widetilde{p}}^{\,intra}$ are the intra-zonal losses that are parameters calculated offline. In case that the losses are not considered implicit for certain interconnectors (e.g. the AC or the HVDC interconnectors), then the corresponding losses are not elements of $\boldsymbol{p}^{loss\textsc{ac}}$ or $\boldsymbol{p}^{loss\textsc{dc}}$, but are included in $\boldsymbol{\widetilde{p}}^{\,intra}$.

Finally, nodal (and zonal) prices are calculated as:
\begin{equation}\label{eq:4_4_4}
    \begin{split}
        \text{LMP}_n = \lambda + & \sum_l \text{PTDF}_{n,l} (\underline{\mu}^{\textsc{AC}}_{\,l} - \overline{\mu}^{\,\textsc{AC}}_l)\,\,+ \\
        & \sum_l \alpha^{\textsc{ac}}_l \,\text{PTDF}_{n,l} ( \sigma^{\textsc{ac},-}_l - \sigma^{\textsc{ac},+}_l)
    \end{split}
\end{equation}
\vspace{-0.6em}

As problem \eqref{eq:3_1_1} aims at minimizing total generation costs, and losses are considered in the power balance equation \eqref{3_1:balance}, the optimization will try to minimize losses. This will lead one of the two inequalities \eqref{eq:4_4_1} to become binding, and, as a result, accurately represent \eqref{eq:4_3_2} or \eqref{eq:4_3_3}. However, in case of negative LMPs, the solver might decide to create artificial losses in order to reduce the system cost \cite{5_12, 5_13, 5_14}. In order to address this issue, in the following the causes of negative LMPs are discussed and a condition for the relaxation to be exact is provided.

\vspace{0.3em}
\noindent
\textbf{Lemma} \textit{For a DC optimal power flow problem as \eqref{eq:3_1_1}, negative LMPs can occur when the cost functions of generators are negative (negative bids) \cite{5_12, 5_13, 5_14, 5_15}, when the total demand is less than the total minimum generation \cite{5_12}, when inter-temporal constraints are included \cite{5_15} and, in case of congestion, when the difference between the marginal costs of production of the marginal generators is big enough and some of the node injections contribute to relieve the congestion.}

\vspace{0.3em}
Due to the equal distribution of losses between the two connected nodes (respectively $l(f)$ and $l(t)$), artificial losses are created when the average price of the two nodes is negative. In the following, we prove that if average prices are always positive, no artificial losses are created.

\vspace{0.3em}
\noindent
\textbf{Proposition} \textit{If the original problem is feasible and  $\frac{1}{2}(\text{LMP}_{l(f)}+\text{LMP}_{l(t)}) > 0$, $\forall l$, then the inclusion of loss functions in the form of two inequalities does not create artificial losses in the system, and the two inequality constraints represent in an exact way the linearized loss functions.}

\vspace{0.3em}
\noindent
\textit{Proof}  From the stationarity conditions of the problem we have:
\begin{equation}\label{eq:4_4_5}
    \frac{1}{2}(\text{LMP}_{l(f)}+\text{LMP}_{l(t)})-(\sigma^+_l+\sigma^-_l) = 0,\quad \forall l
\end{equation}
Being $\sigma^+_l$ and $\sigma^-_l$ Lagrangian multipliers associated with inequality constraints, they are always non-negative. Moreover, only one of the two inequality constraints can be binding at a time, depending on the direction of the flow. For this reason, given that the average price between the two connected node is greater than zero, the stationarity condition is satisfied only if either $\sigma^+_l$ or $\sigma^-_l$ are greater than zero. This means that one of the two inequality constraints is binding, ensuring that no artificial losses are created. $\Box$ 
\noindent

Negative prices are occasionally seen in different markets \cite{5_16, 5_17}. In case of negative prices, many market clearing software use a \textit{Branch and Bound algorithm} to limit losses to their physical value \cite{5_19}. This is also done in case of loop flows on parallel cables \cite{2_9, 5_19}. Another option is the \textit{Big M method}, introducing a binary variable per line and two continuous variables for the flow in the two directions. Although simple, this method slows down the clearing process. However, we never experienced negative prices in our simulations, and no artificial losses were created. For this reason we used the formulation in \eqref{eq:3_1_1} without any of the proposed methods.

Reconnecting with the analysis in \cite{2_15}, including HVDC loss factors for only HVDC interconnectors might be sub-optimal. Indeed, only by including losses on both AC and HVDC lines, the power flows are distributed in a way that minimizes total losses.

\vspace{0.3em}
\noindent
\textbf{Proposition:} \textit{If AC and HVDC loss factors are included for all transmission lines in the market clearing algorithm, the total losses are minimized and the social welfare is always greater than or equal to the case where no losses or only HVDC (or AC)} losses are considered.
\vspace{0.3em}

\textit{Proof}  Let's call \textit{Problem 1} the optimization problem \eqref{eq:3_1_1} with no loss factors and constant losses $\boldsymbol{\widetilde{p}}^{\,loss\textsc{n}}$, \textit{Problem 2} the optimization problem with only HVDC (or AC) loss factors, i.e. only $\boldsymbol{p}^{loss\textsc{dc}}$ are variables in \eqref{eq:4_4_3} and $\boldsymbol{p}^{loss\textsc{ac}}$ are still parameters (or vice versa), and \textit{Problem 3} the optimization problem with both AC and HVDC loss factors, i.e. all the elements of \eqref{eq:4_4_3} are variables and only $\boldsymbol{\widetilde{p}}^{\,intra}$ is a parameter.

For \textit{Problem 3}, the vector of decision variables is $\boldsymbol{\text{x}_3}=[\boldsymbol{g};\boldsymbol{d};\boldsymbol{f}^{\textsc{dc}};\boldsymbol{p}^{loss\textsc{ac}};\boldsymbol{p}^{loss\textsc{dc}}]$ and the feasible space $\Gamma_3$ is the set of solutions that satisfy Eq.  \eqref{3_1:gen}-\eqref{3_1:balance}, Eq. \eqref{eq:4_4_1} and Eq. \eqref{eq:4_4_4}. Including losses as parameters in \textit{Problem 2} is equivalent to adding a new set of constraints to $\Gamma_3$, fixing $\boldsymbol{p}^{loss\textsc{ac}}$ to a certain value $\boldsymbol{\widetilde{p}}^{\,loss\textsc{ac}}$ calculated offline . This means that $\Gamma_2$ is a subset of $\Gamma_3$. A restriction of the feasible space means that the objective value of \textit{Problem 2} can only be less or equal to the objective value of \textit{Problem 3}.

The question arises what would happen if $\boldsymbol{\widetilde{p}}^{\,loss\textsc{ac}}$ is not feasible for \textit{Problem 3}. In that case, $\boldsymbol{\widetilde{p}}^{\,loss\textsc{ac}}$ is an underestimation of the losses considered in \textit{Problem 3}. Given that \textit{Problem 3} provides a better approximation of the actual losses, the solution of \textit{Problem 2} would require the purchase of additional reserves to cover the losses that were not accounted for in the day-ahead. The cost of such reserves are almost always higher than the day-ahead market. As a result, solution $\boldsymbol{\text{x}^*_3}$ always leads to a higher social welfare and an economic benefit.

Following the same approach, $\boldsymbol{\text{x}^*_1}$ leads to an objective value that is less or equal the objective value obtained with $\boldsymbol{\text{x}^*_3}$, as in \textit{Problem 1} all the losses are fixed to a certain value calculated offline, while in \textit{Problem 3} the total losses are allowed to be minimized. $\Box$ 

\vspace{-1em}
\subsection{Comparison of loss factor formulations}\label{ssec:4_5}

\input{fig_3busex.tex}
\begin{table}[!t]
    \caption{Generator, load and line data}
    \vspace{-0.7em}
    \label{tab:4_bus}
    \centering
    \begin{minipage}[t]{0.24\textwidth}
        \footnotesize
        \centering
        \textsc{generators} \\
        \vspace{0.4em}
        \begin{tabularx}{0.80\textwidth}{|*1{>{\centering\arraybackslash}X} c c|}
            \hline 
            ID & $G^{max}$ & $c$  \Tstrut\\
            & \footnotesize{(MW)} & \footnotesize{(\$/MWh)} \Bstrut\\
            \hline
            $g_1$ & 300 & 20 \Tstrut\\
            $g_2$ & 80 & 10 \Bstrut\\
            \hline
        \end{tabularx}
        
        \vspace{1.2em}
        \textsc{ac lines} \\
        \vspace{0.4em}
        \begin{tabularx}{0.80\textwidth}{|*1{>{\centering\arraybackslash}X} c c|}
        \hline
        Line & $F^{max}_{\textsc{ac}}$ & $B$ \Tstrut\\
        & \footnotesize{(MW)} & \footnotesize{(p.u.)} \Bstrut\\
        \hline 
        1-3 & 200 & 0.106 \Tstrut\Bstrut\\
        \hline
        \end{tabularx}
        \vspace{-0.5em}
    \end{minipage}%
    \begin{minipage}[t]{0.24\textwidth}
        \centering
        \footnotesize
            \textsc{loads}\\
        \vspace{0.4em}
        \begin{tabularx}{0.80\textwidth}{|*2{>{\centering\arraybackslash}X}|}
            \hline 
            ID & $D$ \Tstrut\\
            & \footnotesize{(MW)} \Bstrut\\
            \hline 
            $d$ & 292 \Tstrut\Bstrut\\
            \hline
        \end{tabularx}
        
        \vspace{1.2em}
        \textsc{hvdc lines} \\
        \vspace{0.4em}
        \begin{tabularx}{0.80\textwidth}{|*2{>{\centering\arraybackslash}X}|}
       \hline 
        Line & $F^{max}_{\textsc{dc}}$ \Tstrut\\
        & \footnotesize{(MW)} \Bstrut\\
       \hline
        1-2 & 200 \Tstrut\\
        2-3 & 200 \Bstrut\\
        \hline
    \end{tabularx}
    \end{minipage}
    \normalsize
    \vspace{-0.5em}
\end{table}
\begin{table}[!t]
    \caption{hvdc loss factors}
    \vspace{-0.3em}
    \label{tab:4_LFs}
    \footnotesize
    \centering
    \begin{tabularx}{0.44\textwidth}{|c *1{>{\centering\arraybackslash}X} c|}
        \hline 
        Type & Line 1-2 & Line 2-3  \Tstrut\Bstrut\\
        \hline 
        Constant & $\beta = 0.0348$  & $\beta = 0.0332$ \Tstrut\Bstrut\\
        \multirow{2}{*}{Linear} & $\alpha = 0.0403$ & $\alpha = 0.0373$ \Tstrut\\
        & $\beta = 0.0001$ & $\beta = 0.0010$ \Bstrut\\
        \multirow{6}{*}{PW-linear} & $\alpha_1 = 0.0188$ & $\alpha_1 = 0.0171$ \Tstrut\\
        & $\beta_1 = 0.0095$ & \text{$\beta_1 = 0.0100$ } \\
        & $\alpha_2 = 0.0403$ & \text{$\alpha_2 = 0.0373$ } \\
        & $\beta_2 = -0.0048$ & \text{$\beta_2 = -0.0036$ } \\
        & $\alpha_3 = 0.0618$ & \text{$\alpha_3 = 0.0576$ } \\
        & $\beta_3 = -0.0335$ & \text{$\beta_3 = -0.0306$} \Bstrut\\
        \hline
    \end{tabularx}
    \vspace{-0.9em}
    \normalsize
\end{table}

Consider the three-bus network in \figurename~\ref{fig:4_3bus}. To make this illustrative example general, the term ``bus'' is used to refer to different locations in the network: these might correspond to nodes or to zones depending on the pricing scheme. In addition, the load is considered inelastic. Two different system configurations are analyzed: on the left, generator $g_2$ is located in zone 2 and load $d$ in zone 3, while, on the right, their position is swapped. Generator, load and network data are listed in Table~\ref{tab:4_bus}. To study the differnet properties of the proposed formulations, loss factors are introduced only for HVDC lines (Table~\ref{tab:4_LFs}), as proposed in \cite{2_15}. The base power is $\SI{100}{\mega\watt}$ and the base voltage $\SI{400}{\kilo\volt}$.

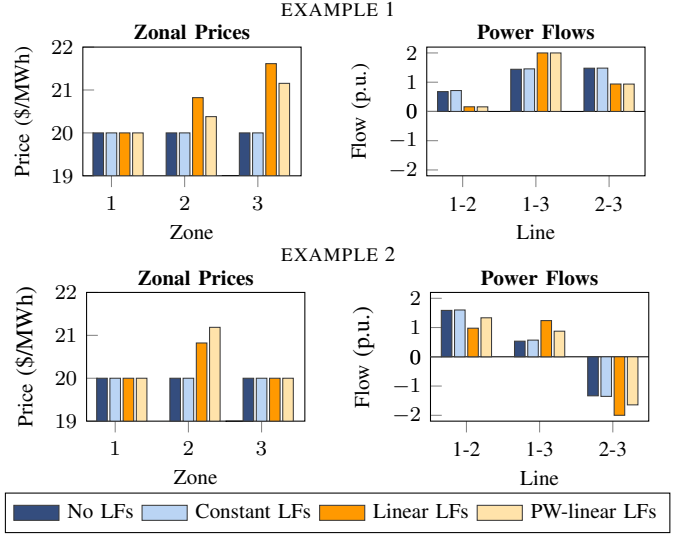
\begin{figure}
    \begin{center}
    \small{\textsc{example }}\footnotesize{1}
    \end{center}
    \vspace{-0.6em}
        \begin{tikzpicture}
            \begin{axis}[%
                width=0.16\textwidth,
                height=0.07\textheight,
                at={(1.20in,0.806in)},
                scale only axis,
                log origin=infty,
                xmin=0.609090909090909,
                xmax=3.59090909090909,
                xtick={1, 2, 3},
                xlabel style={font=\color{white!15!black}},
                xlabel={Zone},
                ymin=19,
                ymax=22,
                ybar=1,
                ylabel style={font=\color{white!15!black}},
                ylabel={Price {(\$/MWh)}},
                ylabel near ticks,
                xtick pos=left,
                ytick pos=left,
                label style={font=\footnotesize},
                every tick label/.append style={font=\footnotesize},
                axis background/.style={fill=white},
                title style={font=\bfseries\footnotesize,yshift=-1.5ex},
                title={Zonal Prices}
                ]
                \addplot[black] table[x index = 0, y index = 1] {
                0.5 0
                3.49 0
                };
                
                \addplot[ybar, bar width=4, fill=blue1, draw=gray1] table[row sep=crcr] {%
                1	20\\
                2	20.0000000000661\\
                3	20.0000000075899\\
                };
                
                \addplot[ybar, bar width=4, fill=blue2, draw=gray1] table[row sep=crcr] {%
                1	20\\
                2	20.0000000000672\\
                3	20.0000000080057\\
                };
                
                \addplot[ybar, bar width=4, fill=orange2, draw=gray1] table[row sep=crcr] {%
                1	20.0000000584777\\
                2	20.821950010623\\
                3	21.6142335424816\\
                };
                
                \addplot[ybar, bar width=4, fill=orange3, draw=gray1] table[row sep=crcr] {%
                1	20.0000000327184\\
                2	20.3793096444664\\
                3	21.1547504536608\\
                };
                
            \end{axis}
    
            \begin{axis}[%
                width=0.16\textwidth,
                height=0.07\textheight,
                at={(3in,0.806in)},
                scale only axis,
                bar shift auto,
                xmin=0.509090909090909,
                xmax=3.49090909090909,
                xtick={1,2,3},
                xticklabels={{1-2},{1-3},{2-3}},
                xlabel style={font=\small\color{white!15!black}},
                xlabel={Line},
                extra y ticks = 0,
                extra y tick style ={ grid=major,major grid style={draw=black}}, 
                ymin=-2.20,
                ymax=2.20,
                ylabel style={font=\color{white!15!black}},
                ylabel={Flow {(p.u.)}},
                ylabel near ticks,
                ybar=1,
                xtick pos=left,
                ytick pos=left,
                label style={font=\footnotesize},
                every tick label/.append style={font=\footnotesize},
                axis background/.style={fill=white},
                title style={font=\bfseries\footnotesize,yshift=-1.5ex},
                title={Power Flows},
                ]
                
                \addplot[ybar, bar width=4, fill=blue1, draw=gray1] table[row sep=crcr] {%
                1	0.678314925936213\\
                2	1.44168506819167\\
                3	1.47831493180833\\
                };
                
                \addplot[ybar, bar width=4, fill=blue2, draw=gray1] table[row sep=crcr] {%
                1	0.715230833952818\\
                2	1.4553923604154\\
                3	1.4812081395846\\
                };
                
                \addplot[ybar, bar width=4, fill=orange2, draw=gray1] table[row sep=crcr] {%
                1	0.159292211065301\\
                2	1.99999998041561\\
                3	0.938029050537985\\
                };
            
                \addplot[ybar, bar width=4, fill=orange3, draw=gray1] table[row sep=crcr] {%
                1	0.15760709079857\\
                2	2.00000000737394\\
                3	0.935687047491055\\
                };
                
            \end{axis}
    \end{tikzpicture}
    \vspace{-1em}        
    \begin{center}
    \small{\textsc{example }}\footnotesize{2}
    \end{center}
    \vspace{-0.6em} 
        \begin{tikzpicture}
            \begin{axis}[%
                width=0.16\textwidth,
                height=0.07\textheight,
                at={(1.20in,0.806in)},
                scale only axis,
                log origin=infty,
                xmin=0.609090909090909,
                xmax=3.59090909090909,
                xtick={1, 2, 3},
                xlabel style={font=\color{white!15!black}},
                xlabel={Zone},
                ymin=19,
                ymax=22,
                ybar=1,
                ylabel style={font=\color{white!15!black}},
                ylabel={Price \small{(\$/MWh)}},
                ylabel near ticks,
                xtick pos=left,
                ytick pos=left,
                label style={font=\footnotesize},
                every tick label/.append style={font=\footnotesize},
                axis background/.style={fill=white},
                title style={font=\bfseries\footnotesize,yshift=-1.5ex},
                title={Zonal Prices}
                ]
                \addplot[black] table[x index = 0, y index = 1] {
                0.5 0
                3.49 0
                };
                
                \addplot[ybar, bar width=4, fill=blue1, draw=gray1] table[row sep=crcr] {%
                1	20\\
                2	20.0000000074746\\
                3	20.0000000071716\\
                };
                
                \addplot[ybar, bar width=4, fill=blue2, draw=gray1] table[row sep=crcr] {%
                1	20\\
                2	20.0000000087961\\
                3	20.0000000084817\\
                };
                
                \addplot[ybar, bar width=4, fill=orange2, draw=gray1] table[row sep=crcr] {%
                1	20.0000000949772\\
                2	20.8219500413996\\
                3	20.0000002122031\\
                };
                
                \addplot[ybar, bar width=4, fill=orange3, draw=gray1] table[row sep=crcr] {%
                1	20.0000002270186\\
                2	21.1864720954167\\
                3	20.000000317239\\
                };
                
            \end{axis}
    
            \begin{axis}[%
                width=0.16\textwidth,
                height=0.07\textheight,
                at={(3in,0.806in)},
                scale only axis,
                bar shift auto,
                xmin=0.509090909090909,
                xmax=3.49090909090909,
                xtick={1,2,3},
                xticklabels={{1-2},{1-3},{2-3}},
                xlabel style={font=\color{white!15!black}},
                xlabel={Line},
                extra y ticks = 0,
                extra y tick style ={ grid=major,major grid style={draw=black}}, 
                ymin=-2.20,
                ymax=2.20,
                ylabel style={font=\color{white!15!black}},
                ylabel={Flow \small{(p.u.)}},
                ylabel near ticks,
                ybar=1,
                xtick pos=left,
                ytick pos=left,
                label style={font=\footnotesize},
                every tick label/.append style={font=\footnotesize},
                axis background/.style={fill=white},
                title style={font=\bfseries\footnotesize,yshift=-1.5ex},
                title={Power Flows},
                legend columns=4,
                legend style={at={(-1.95,-0.87)}, anchor=south west, legend cell align=left, align=left, draw=black, font=\footnotesize},
                every axis legend/.append style={column sep=0.1em}
                ]
                
                \addplot[ybar, bar width=4, fill=blue1, draw=gray1, area legend] table[row sep=crcr] {%
                1	1.58607439719173\\
                2	0.533925600324083\\
                3	-1.33392560280827\\
                };
                \addlegendentry{No LFs}
                
                \addplot[ybar, bar width=4, fill=blue2, draw=gray1, area legend] table[row sep=crcr] {%
                1	1.60046270355138\\
                2	0.570160493750203\\
                3	-1.35355999644862\\
                };
                \addlegendentry{Constant LFs}
                
                \addplot[ybar, bar width=4, fill=orange2, draw=gray1, area legend] table[row sep=crcr] {%
                1	0.977566342134092\\
                2	1.23785577256744\\
                3	-1.99999975422931\\
                };
                \addlegendentry{Linear LFs}
            
                \addplot[ybar, bar width=4, fill=orange3, draw=gray1, area legend] table[row sep=crcr] {%
                1	1.33333342667821\\
                2	0.87517635504317\\
                3	-1.64314055914736\\
                };
                \addlegendentry{PW-linear LFs}
                
            \end{axis}
        \end{tikzpicture}
    \caption{Zonal prices and line flows.}
    \label{fig:4_LMP}
    \vspace{-1em}
\end{figure}

To compare the impact of the different loss factor formulations, the optimization problem \eqref{eq:3_1_1} is solved four times. The first time, no HVDC loss factors are included. The other times, constraint \eqref{eq:4_4_3} is included together with, respectively, constant, linear and piecewise linear loss factors. \figurename~\ref{fig:4_LMP} shows the different prices and power flows obtained with the four formulations.

In Example 1, most of the power flows from bus 1 to bus 3. With this configuration, the power has two possible paths, either over one AC interconnector or over two HVDC lines. When the market is cleared without loss factors, no distinction is made between HVDC and AC lines and, thus, there are several power flow solutions for the market equilibrium. If constant loss factors are introduced, losses still do not depend on the power flows, so prices and flows remain unchanged. 

\onehalfspacing
The situation changes when linear and piecewise linear loss factors are introduced. Indeed, losses are now a function of the power flow, and thus the more the HVDC lines are used, the higher the losses are. For this reason, the use of HVDC lines is limited to the cases when the AC capacity constraint violation cannot be resolved by any other measure. In addition, a price difference is forced between buses 1-2 and buses 2-3 when the HVDC line is used. These price differences are functions of the linear coefficients of losses and can be calculated as:
\vspace{0.3em}
\small
\begin{equation}\label{eq:4_9}
    \text{LMP}_2 = \frac{1+0.5\,\alpha^{\textsc{dc}}_{1}}{1-0.5\,\alpha^{\textsc{dc}}_{1}} \cdot \text{LMP}_1 \quad \text{LMP}_3 = \frac{1+0.5\,\alpha^{\textsc{dc}}_{2}}{1-0.5\,\alpha^{\textsc{dc}}_{2}} \cdot \text{LMP}_2.
\end{equation}
\vspace{-2.5em}
\normalsize \singlespacing \noindent %
These equations are derived from the KKT optimality conditions, and give the relation between the lagrangian multipliers associated with the power balance equations. Once the limit of line 1-3 is reached, the only way to supply the load is through the two HVDC lines. An increase of consumption $\Delta d$ at bus 3 would correspond to an increase of generation equal to $\Delta d$ plus the losses, and thus it would be more expensive than an equal increase at bus 1 or 2. In case of piecewise linear loss factors, the coefficients appearing in \eqref{eq:4_9} are the linear coefficients of the binding loss functions. It should be mentioned that \eqref{eq:4_9} depends on the direction of the HVDC flows. In case of opposite flow between e.g. zone 1 and 2, then the signs in \eqref{eq:4_9} will be opposite, i.e. $\text{LMP}_2 = [(1-0.5\,\alpha^{\textsc{dc}}_{1})/(1+0.5\,\alpha^{\textsc{dc}}_{1})]\text{LMP}_1$.

In Example 2, the load is moved to bus 2 and $g_2$ to bus 3. Now both paths for supplying the load include an HVDC line. As \figurename~\ref{fig:4_LMP} shows, with no loss factors or with constant loss factors the market outcome is very similar: in the first case no distinction between AC and HVDC lines is made, and in the second case losses do not depend on the flows and thus prices and flows remain unchanged. Again, the situation is different when linear or piecewise linear loss factors are introduced. With linear loss factors, the slope of the loss function determines the path that results in less losses. Indeed, the power flow over line 2-3 is equal to its capacity, while only the remaining power is supplied through line 1-2. With piecewise linear loss factors, the slope of the loss function changes depending on the flow. For this reason, the solver identifies the least costly path by moving back and forth from one line to the other, when the slope of the loss function changes. In this way the two lines are used in a more efficient way, and the price difference, although greater, reflects better the cost of losses.
\vspace{-0.5em}
\section{Intra-zonal losses}\label{sec:5}
\begin{figure}[!t]
\centering
\includegraphics[trim = 1.9cm 20cm 11cm 2cm,clip,width=0.48\textwidth]{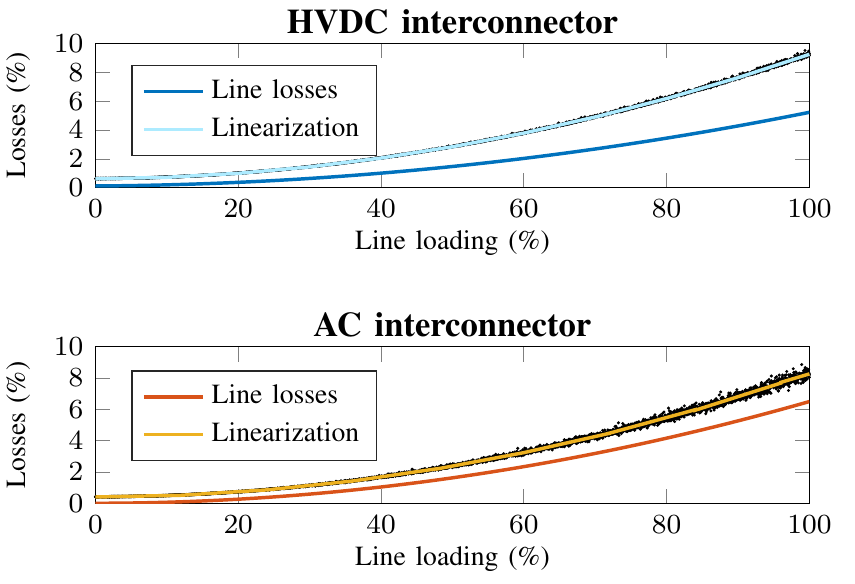}%
\vspace{-1em}
\caption{Inclusion of intra-zonal losses on HVDC and AC loss factors. Line losses are the losses on the interconnector, black dots are the intra-zonal losses due to cross-border flows. The linearization is made with 10 segments.}
\label{fig:5_LFs}
\vspace*{-1.3em}
\end{figure}

In our initial investigations, and as we also demonstrate in our numerical tests in Section~\ref{sec:6}, the social welfare does not always increase by introducing AC and HVDC loss factor in zonal pricing markets. This is due to the fact that intra-zonal losses are not variables in the optimization problem, and thus cannot be minimized. This chapter introduces a method for considering a part of intra-zonal losses, which is produced by cross-border flows, in the calculation of loss factors. 

Intra-zone losses are caused by both cross-border and internal flows. Loss factors are meant to account only for the losses due to inter-zonal flows, so the internal losses have to be excluded from the calculation. We calculated the new loss factors as follows. Two zones connected are considered at a time: in Zone 1 all the loads are removed, while in Zone 2 all generators are removed. The statistical population of losses is calculated running 10'000 AC power flows, where different generation patterns and load conditions are considered. The same procedure is then repeated inverting generation and consumption in the two areas. We repeat these two steps for all the zones connected by AC lines. Ideally, the losses between any two zones would have been the result of the superposition of the losses found for each pair of zones. However, in that case we account for the losses in each transit zone more than once. As a result, we carry out a similar analysis considering the whole system, and estimate a correction factor that we introduce in order to avoid accounting for the same losses twice. Finally, we use the least-squares methods to linearize the losses.

\vspace{-0.7em}
\section{Case study}\label{sec:6}
\input{fig_4area.tex}

We compare four market outcomes (with no loss factors, with only HVDC loss factors, with only AC loss factors and with AC and HVDC loss factors) on a modified version of the IEEE 3-area RTS `96 Test Case \cite{7_1} using YALMIP \cite{7_2} and MOSEK \cite{7_3}. The numerical tests are carried out considering both nodal and zonal pricing mechanisms. The market is cleared for every hour of a year, resulting in a time series of 8760 market outcomes. The test case is modified as described in \cite{7_4}, including a fourth zone, three HVDC interconnectors and some wind farms. To modify energy prices over the year, load consumption and wind generation vary in each snapshot. This leads to different import-export situations for each zone, resulting in a more realistic analysis.

\begin{figure}[!t]
\centering
\includegraphics[trim = 1.7cm 20.7cm 11cm 1.8cm,clip,width=0.48\textwidth]{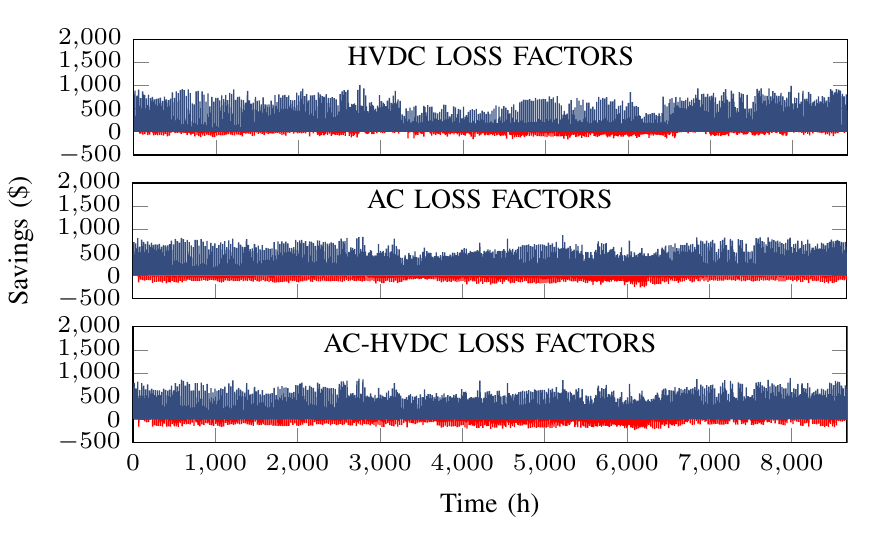}%
\vspace*{-1em}
\caption{FBMC: comparison of market outcomes. Loss factors calculated based on line losses and linear approximation.}
\label{fig:6_TSL}
\vspace*{-1em}
\end{figure}

\begin{figure}[!t]
\centering
\includegraphics[trim = 1.7cm 14.5cm 11cm 8cm,clip,width=0.48\textwidth]{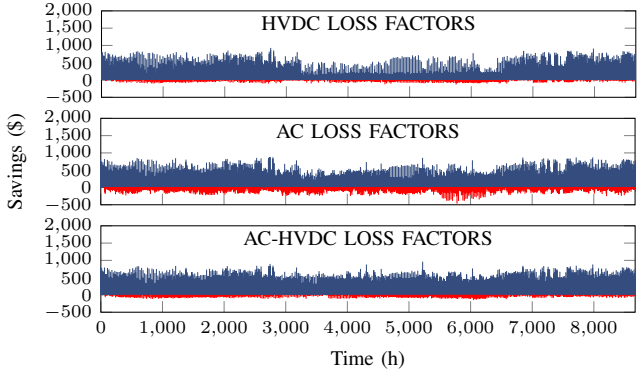}%
\vspace*{-1em}
\caption{FBMC: comparison of market outcomes. Loss factors calculated based on line losses and piecewise linear approximation.}
\label{fig:6_TSPWL}
\vspace*{-1.2em}
\end{figure}

\begin{figure}[!t]
\centering
\includegraphics[trim = 11cm 20.7cm 1.7cm 1.8cm,clip,width=0.48\textwidth]{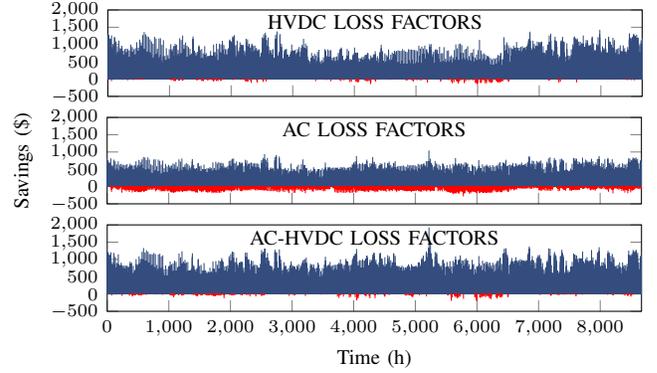}%
\vspace*{-1em}
\caption{FBMC: comparison of market outcomes. Loss factors calculated based on line losses, intra-zonal losses and  piecewise linear approximation.}
\label{fig:6_TSCB}
\vspace*{-1em}
\end{figure}

\begin{figure}[!t]
\centering
\includegraphics[trim = 11cm 14.6cm 1.7cm 7.9cm,clip,width=0.48\textwidth]{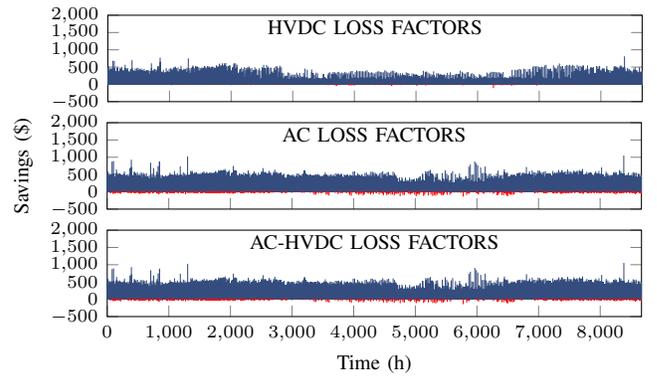}%
\vspace*{-1em}
\caption{Nodal pricing: comparison of market outcomes. Loss factors calculated based on piecewise linear approximation.}
\label{fig:6_TSN}
\vspace*{-1.2em}
\end{figure}

To investigate the impact of the introduction of loss factors on the social welfare, four simulations are run. The first three times, a zonal pricing scheme is considered and the 96-bus system is reduced to a 4-area system for FBMC, as shown in \figurename~\ref{fig:6_4areasys}. The fourth time, a nodal pricing scheme is considered and the whole network is included in the model. In all the simulations, the following procedure is applied:
\begin{enumerate}
\item Losses are calculated offline. The market is cleared without considering losses. Subsequently, the determined generation and load setpoints are used to calculate the losses;
\item Four optimization problems are solved. In the first, losses are considered as constant parameters. In the second, only HVDC losses are variables in the optimization problem, while all AC losses are parameters. In the third, only AC losses are variables in the optimization problem, while all HVDC losses are parameters. In the fourth, the losses in all the interconnectors are variables and only the intra-zonal losses are still calculated offline;
\item The social welfare is computed according to the three market outcomes. The differences are plotted in \figurename~\ref{fig:6_TSL}, \ref{fig:6_TSPWL}, \ref{fig:6_TSCB} and \ref{fig:6_TSN}. 
\end{enumerate}

For the first simulation, loss factors are calculated considering only losses on the interconnectors. The quadratic loss functions are linearized with linear loss factors considering two fixed points (no-flow and 60\% line loading), similarly to the proposal of Nordic TSOs in \cite{2_15}. No intra-zonal losses are considered in the loss factor calculation. \figurename~\ref{fig:6_TSL} shows the differences in social welfare between the base case (no loss factors) and the other optimization problems, with implicit grid loss implemented respectively on only HVDC, only AC and both AC and HVDC interconnectors.


Overall, the inclusion of losses in the market clearing has a positive impact. However, in many simulations, the social welfare is decreased. First of all, the estimation of losses is not very accurate. As explained, the market is cleared without losses, and these are calculated based on the set points of generators and loads. Including losses as constant loads alters the flows, thus the losses calculated are not precise. This happens also in reality, since losses are forecast. The underestimation of losses requires the purchase of this power in the balancing market resulting in higher costs, thus having a representation of losses in the market clearing is a good solution. Moreover, internal losses are not considered when defining the power exchanges: different power flows on the interconnectors might result in higher amount of internal losses, and thus in higher costs.

For the second simulation, still no intra-zonal losses are considered, but line losses are now linearized using piecewise linear functions with segments of the length of 60 MW. By doing so, flows are distributed in a better way, since the solver moves back and forth from one loss function to the other determining the least costly path. Moreover, the quadratic loss functions are approximated in a better way. This leads to a partial improvement of the results, as shown in \figurename~\ref{fig:6_TSPWL}. However, with only AC loss factors, the situation seems to get worse. This is because, in most of the situations, Zone 1 and 4 are exporting. With AC loss factors, the power from Area 4 is rerouted through Area 1, creating more losses on the interconnectors and causing also more internal losses, while, with DC loss factors, HVDC lines are the only path from Zone 1 to the others, thus HVDC lines are used anyways, and with both loss factors, all the lines are used in an optimal way.

In the third simulation, additionally to line losses, intra-zonal losses due to cross-border flows are included in the calculation of loss factors. Different points of connection might result in different internal losses, and this is now considered. As a consequence, the results are further improved, as shown in figure \figurename~\ref{fig:6_TSCB}. 


Finally, a similar analysis is carried out with a nodal pricing scheme, considering the whole network. Contrary to the other simulations, there are now 156 AC lines and 3 HVDC lines. This means that the possibility of controlling the flows is reduced, and thus the amount of savings as well. However, the representation of losses in the market clearing is more accurate, and there are fewer situations where the social welfare is decreased. Also, one would expect to have only positive increments with AC and HVDC loss factors, however the inaccuracy of the loss estimation procedure affects the results.

\tablename~\ref{tab:6_comparison} presents the total savings in the different situations. In all simulations, the inclusion of both AC and HVDC loss factors gives the greatest benefit, as theoretically expected.

\begin{table}[!t]
    \vspace{-0.5em}
    \caption{total savings (M\$)}
    \vspace{-0.4em}
    \label{tab:6_comparison}
    \small
    \centering
    \begin{tabular}{|cccc|}
        \hline 
        & HVDC LF & AC LF & AC-HVDC LF \Tstrut\\
        \hline 
        Simulation 1 & 1.16    & 1.12  & 1.38 \Tstrut\\
        Simulation 2 & 1.05    & 0.90  & \text{$1.51$ } \\
        Simulation 3 & 2.52    & 1.23  & \text{$3.10$ } \\
        Simulation 4 & 0.99    & 1.77  & 1.81 \Bstrut\\
        \hline
    \end{tabular}
    \normalsize
\end{table}
\vspace{-0.5em}

\section{Conclusion}\label{sec:7}
The introduction of loss factors for HVDC lines, also called implicit grid loss calculation, has been proposed by the TSOs of Nordic Capacity Calculation Region to avoid HVDC flows between zones with zero price difference. Currently, it is under investigation for real implementation in the market clearing algorithm. In this paper, we have introduced a rigorous framework to assess the impact of the shift towards implicit grid losses, considering the introduction of loss factors for different interconnectors. We develop different loss factor formulations and study their main properties on a representative test system. We find that although the introduction of HVDC loss factors is in general positive, it may lead to a decrease of the social welfare for a non-negligible amount of time as it disproportionately increases the AC losses. For zonal pricing markets, this might happen also when implicit grid losses are implemented in all interconnectors because of intra-zonal losses. To counter that, we introduce a methodology to estimate loss factors based on statistical analysis and linear regression. We prove theoretically that the introduction of both AC and HVDC loss factors in market clearing is guaranteed to increase the social welfare. We confirm our results through numerical tests in a representative test system. 


\bibliographystyle{myIEEEtran.bst}

\vspace{-1em}
\vspace{-2em}
\begin{IEEEbiography}[{\includegraphics[width=1in,height=1.25in,clip]{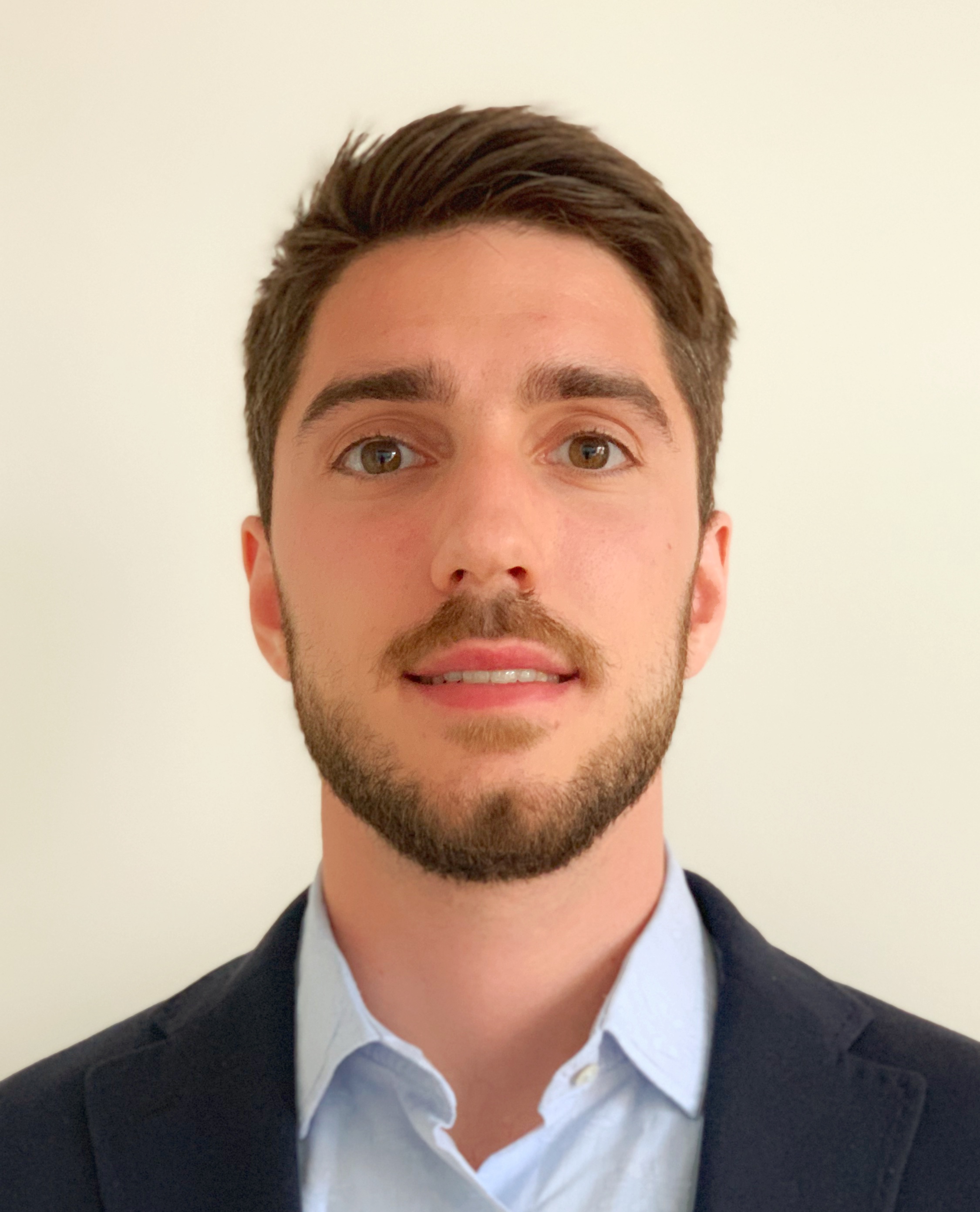}}]{Andrea Tosatto}
received the M.Sc. degree in electrical engineering from the Royal Institute of Technology (KTH), Stockholm, Sweden, and a second M.Sc. degree in electrical engineering from the Institute Polytechnique de Grenoble, Grenoble, France. He is working toward the Ph.D. degree at Centre for Electric Power and Energy, Department of Electrical Engineering, Technical University of Denmark (DTU). His research interests include convex optimization in power systems, applied game theory, and market integration of multiarea AC/HVDC grids.
\end{IEEEbiography}

\begin{IEEEbiography}[{\includegraphics[width=1in,height=1.38in,clip,keepaspectratio]{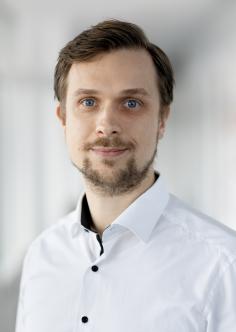}}]{Tilman Weckesser}
received the M.Sc. and Ph.D. degrees from the Technical University of Denmark (DTU), Lyngby, Denmark, in 2011 and 2015, respectively. He was with the University of Lige, Belgium, as a Postdoctoral Researcher and with the Technical University of Denmark as an Assistant Professor. Currently, he is Consultant with the Danish Energy Association in the Grid Technology Department. His research interests include electric power system dynamics, stability, and electric mobility.
\end{IEEEbiography}
\vspace{-3em}
\begin{IEEEbiography}[{\includegraphics[width=1in,height=1.25in,clip]{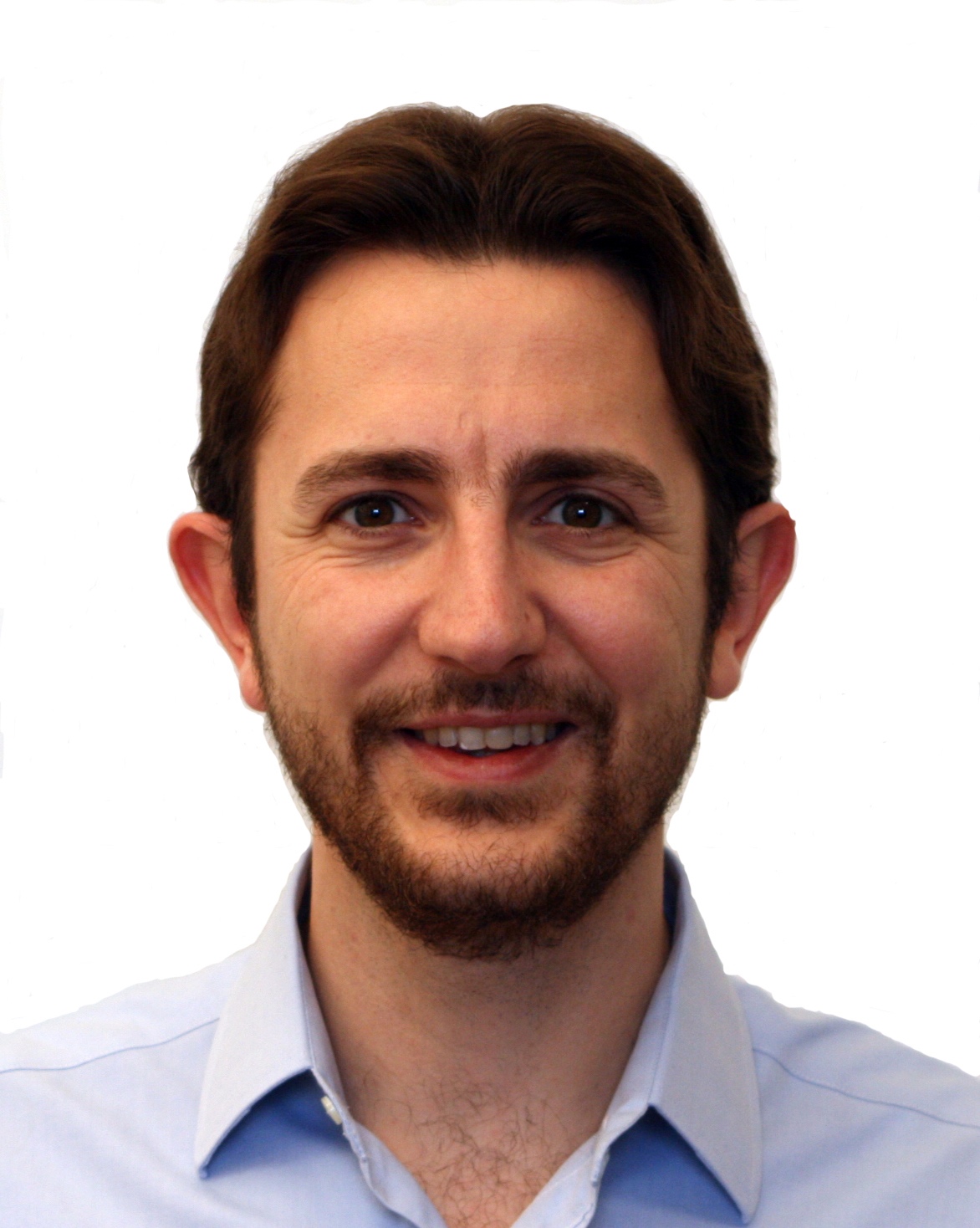}}]{Spyros Chatzivasileiadis}
(S'04, M'14) received the Diploma in electrical and computer engineering from the National Technical University of Athens (NTUA), Greece, in 2007, and the Ph.D. degree from ETH Zurich, Switzerland in 2013. He was a Postdoctoral Researcher with the Massachusetts Institute of Technology (MIT), USA and with Lawrence Berkeley National Laboratory, USA. In March 2016 he joined
the Center for Electric Power and Energy, Technical University of Denmark (DTU). He is an Associate Professor with DTU. He is currently working on power system optimization and control of AC and HVDC grids, and data-driven optimization and security assessment of transmission and distribution grids.
\end{IEEEbiography}

\balance

\end{document}